\documentclass[review]{elsarticle}

\usepackage{lineno,hyperref, amsmath, amssymb, graphics, arydshln, amsthm, enumitem, url}
\newtheorem{theorem}{Theorem}
\newtheorem{remark}{Remark}

\modulolinenumbers[5]

\begin{document}

\begin{frontmatter}

\title{Narayana, Mandelbrot, and \\ A New Kind of Companion Matrix}

\author{Eunice Y. S. Chan}
\ead{echan295@uwo.ca}
\author{Robert M. Corless}
\ead{rcorless@uwo.ca}
\address{Ontario Research Centre for Computer Algebra \\ Department of Applied Mathematics \\ Western University \\ London, Ontario}

\begin{abstract}
We demonstrate a new kind of companion matrix, for polynomials of the form $c(\lambda) = \lambda a(\lambda)b(\lambda) + c_0$ where upper Hessenberg companions are known for the polynomials $a(\lambda)$ and $b(\lambda)$. This construction can generate companion matrices with smaller entries than the Fiedler or Frobenius forms. This generalizes Piers Lawrence's Mandelbrot companion matrix. We motivate the construction by use of Narayana-Mandelbrot polynomials, which are also new to this paper.
\end{abstract}


\end{frontmatter}


\section{Introduction}

Sequence A000930 of the Online Encyclopedia of Integer Sequences, Narayana's cows sequence, begins
\begin{equation}
	1, 1, 1, 2, 3, 4, 6, 9, 13, 19, \ldots
\end{equation}
and is generated by $r_n = r_{n-1} + r_{n-3}$ \cite{A000930}. The connection to cows is that an ideal cow produces a calf every year, starting in its fourth year. Narayana was a mathematician in 14th century India. Various facts are known for this sequence, which is similar to the Fibonacci sequence: for instance, the generating function is $1/(1-x-x^3)$. Many references are given in the OEIS, but see also~\cite{sloane1999my}.

Recently, we generalized the Mandelbrot polynomials
\begin{equation}
	p_{n+1} = zp_n^2 + 1 \quad p_0=0
\end{equation}
to the Fibonacci-Mandelbrot polynomials
\begin{equation}
	q_{n+1} = zq_nq_{n-1} + 1 \quad q_0=0, q_1=1
\end{equation}
and generalized Piers Lawrence's supersparse\footnote{A matrix is supersparse if, it is sparse and its nonzero elements are drawn from a small set, e.g. $\{-1, 1\}$} companion matrix for $p_n$ \cite{piers} to an analogous one for $q_n$. See \cite{corless2014graduate}, \cite{chan2016thesis} and \cite{chan2016fibonacci} for details, though we summarize these constructions below.

\begin{figure}
	\centering
	\includegraphics[width = 1 \textwidth]{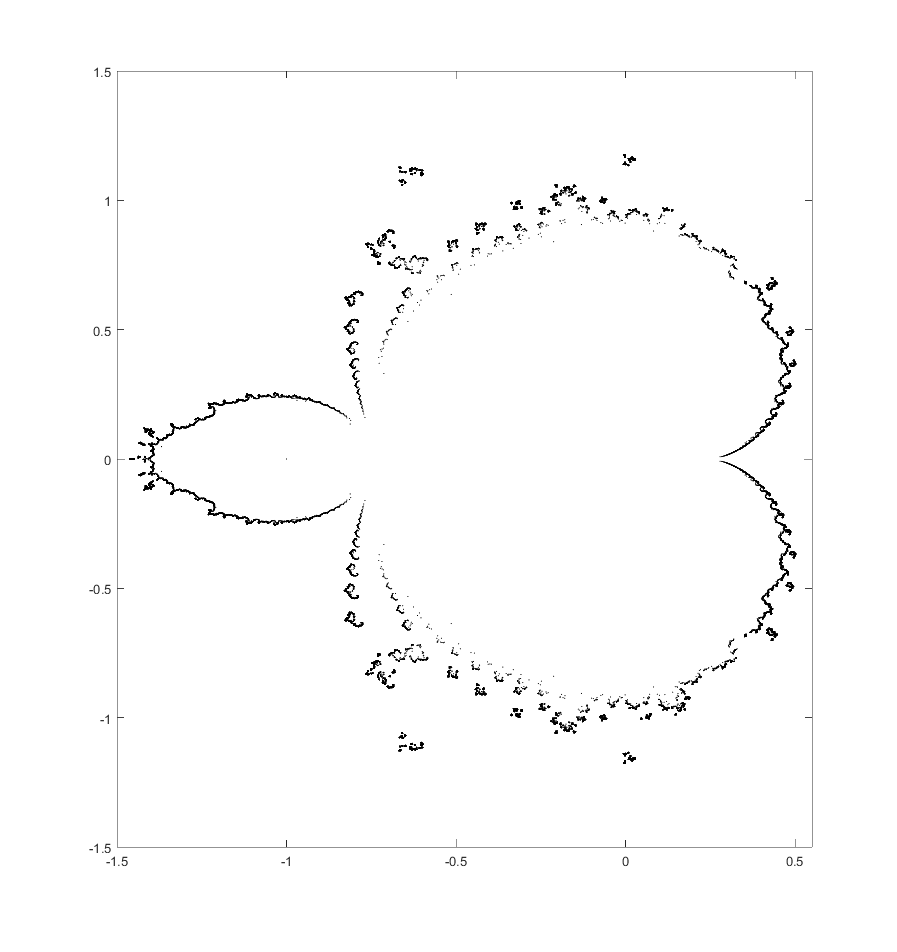}
	\caption{Roots of Narayana-Mandelbrot polynomial, $r_{36}(z)$. The degree of $r_{36}(z)$ is 395032.}
\end{figure}

In this paper, we define the Narayana-Mandelbrot polynomials by $r_0 = 0, r_1 = r_2 = 1$ and
\begin{equation}
	r_{n+1} = zr_nr_{n-2} + 1
\end{equation}
for $n \geq 2$. We give some basic facts about these polynomials, and we construct a recursive family of companion matrices $\mathbf{R}_n$, i.e. such that
\begin{equation}
	r_n(z) = \det(z\mathbf{I} - \mathbf{R}_n).
\end{equation}
The $\mathbf{R}_n$ will be seen to be supersparse. We prove that the construction is valid by using induction and the Schur determinantal formula.

The surprising analogy between all three families of supersparse companions led us to conjecture and prove the following.
\begin{theorem}
	\sl Suppose $a(z) = \det(z\mathbf{I} - \mathbf{A})$, $b(z) = \det(z\mathbf{I} - \mathbf{B})$, and both $\mathbf{A}$ and $\mathbf{B}$ are upper Hessenberg matrices with nonzero subdiagonal entries, and
	\begin{equation}
		\alpha = \cfrac{1}{\left(\prod^{d_a-1}_{j=1}a_{j+1,j}\right)\left(\prod^{d_b-1}_{j=1}b_{j+1,j}\right)}
	\end{equation}
is the reciprocal of the product of the subdiagonal entries of $\mathbf{A}$ and $\mathbf{B}$, and $d_a = \deg_za$ and $d_b=\deg_zb$, so the dimension of $\mathbf{A}$ is $d_a \times d_a$ and the dimension of $\mathbf{B}$ is $d_b \times d_b$. Suppose both $d_a$ and $d_b$ are at least 1. Then if
\begin{equation}
	\mathbf{C} = \left[
	\begin{array}{ccc}
		\mathbf{A} & & -\alpha c_0 \mathbf{c_a}\mathbf{r_b} \\
		-\mathbf{r_a} & 0 & \\
		& -\mathbf{c_b} & \mathbf{B}
	\end{array}
	\right]
\end{equation}
where $\mathbf{r_a} = \left[ \begin{array}{cccc}0 & 0 & \cdots & 1 \end{array}\right]$ of length $d_a$, $\mathbf{c_b} = \left[\begin{array}{cccc} 1 & 0 & \cdots & 0\end{array}\right]^T$ of length $d_b$, we have
\begin{equation}
	c(z) = \det\left(z\mathbf{I} - \mathbf{C}\right) = z\cdot a(z)b(z) + c_0.
\end{equation}
\end{theorem}

\begin{remark}
	\sl Proving this theorem automatically proves the validity of the constructions of the supersparse companion matrices for $p_n$, $q_n$, and $r_n$.
\end{remark}

\begin{remark}
	\sl Starting with a polynomial $c(z)$, we see that there are potentially many such $a(z)$ and $b(z)$. This freedom may be quite valuable or, it may be an obstacle.
\end{remark}

\begin{proof}
	Partition
	\begin{equation}
		z\mathbf{I} - \mathbf{C} = \left[ \begin{array}{c:c} \mathbf{C}_{11} & \mathbf{C}_{12}\\ \hdashline \mathbf{C}_{21} & \mathbf{C}_{22} \end{array} \right]
	\end{equation}
	where $\mathbf{C}_{22} = z\mathbf{I} - \mathbf{B}$ is nonsingular if $z$ is not an eigenvalue of $\mathbf{B}$, i.e. $b(z) \neq 0$. Later we will remove this restriction. Also,
	\begin{equation}
		\mathbf{C}_{21} = \left[ \begin{array}{ccc} & & 1 \\ & &  \\& & \end{array} \right]
	\end{equation}
	is $d_b \times (d_a + 1)$ and has only one nonzero element, which is a $1$ in the upper right corner. Next,
	\begin{equation}
		\mathbf{C}_{12} = \left[ \begin{array}{ccc} & &\alpha c_0 \\ & & \\ & & \end{array}\right]
	\end{equation}
	is $(1 + d_a) \times d_b$ and again has only one nonzero element, $\alpha c_0$ in the upper right corner. [In fact, $c_0$ can be zero.] This leaves
	\begin{equation}
		\newcommand*{\tempb}{\multicolumn{1}{:c}{\begin{array}{c}0 \\ \vdots \\ 0 \\ 0\end{array}}}
		\mathbf{C}_{11}=\left[
		\begin{array}{cc}
			z\mathbf{I} - \mathbf{A} &\tempb \\ \cdashline{1-1}
			\begin{array}{cccc}
				& & &1
			\end{array}
			& z
		\end{array}\right]
	\end{equation}
	which is $d_a + 1$ by $d_a + 1$.
	
	The Schur factoring is
	\begin{equation}
		\left[
		\begin{array}{cc}
			\mathbf{C}_{11} & \mathbf{C}_{12} \\
			\mathbf{C}_{21} & \mathbf{C}_{22}
		\end{array}
		\right]
		=
		\left[
		\begin{array}{cc}
			\mathbf{I} & \mathbf{C}_{12} \\
			0 & \mathbf{C}_{22}
		\end{array}
		\right]
		\left[
		\begin{array}{cc}
			\mathbf{C}_{11}-\mathbf{C}_{12}\mathbf{C}_{22}^{-1}\mathbf{C}_{21} & 0 \\
			\mathbf{C}_{22}^{-1} & \mathbf{I}
		\end{array}
		\right]
	\end{equation}
	with the computation of the Schur complement $\mathbf{C}_{11} - \mathbf{C}_{12}\mathbf{C}^{-1}_{22}\mathbf{C}_{21}$ going to do most of the work in the proof. The Schur determinantal formula \cite[Chapter 12]{hogben2014handbook} is then
	\begin{equation}
		\det \mathbf{C} = \det\left(\mathbf{C}_{22}\right)\det\left(\mathbf{C}_{11} - \mathbf{C}_{12}\mathbf{C}_{22}^{-1}\mathbf{C}_{21}\right).
	\end{equation}
	We have the following propositions.
	\begin{enumerate}[start = 0]
	\item $z\mathbf{I} - \mathbf{A}$ and $z\mathbf{I} - \mathbf{B}$ are upper Hessenberg because $\mathbf{A}$ and $\mathbf{B}$ are.
	\item The first $d_a$ columns of $\mathbf{C}_{22}^{-1}\mathbf{C}_{21}$ are zero.
	\item The final column of $\mathbf{C}_{22}^{-1}\mathbf{C}_{21}$ is the solution, say $\vec{v}$, of $\left(z\mathbf{I} - \mathbf{B}\right)\vec{v} = \mathbf{e}_1$. Again, $z\mathbf{I} - \mathbf{B}$ is nonsingular.
	\item By Cramer's rule, the final entry in $\vec{v}$, say $v$, is
	\begin{equation}
		v = \frac{\det\left(\mathbf{C}_{22}\xleftarrow[d_b]{} \mathbf{e}_1\right)}{\det\left(\mathbf{C}_{22}\right)}
	\end{equation}
	where the notation $\mathbf{M} \xleftarrow[k]{} \vec{v}$ means replace the $k$th column of $\mathbf{M}$ with the vector $\vec{v}$ \cite{carlson1992gems}.
	\item Since $\mathbf{C}_{22} = z\mathbf{I} - \mathbf{B}$ is upper Hessenberg,
	\begin{equation}
		\mathbf{C}_{22}\xleftarrow[d_b]{}e_1 = \left[
		\begin{array}{cccccc}
			* & * & * & \cdots & * & 1 \\
			-b_{21} & * & * & \cdots & * & 0 \\
			& -b_{32} & * & & \vdots & \vdots \\
			& & -b_{43} & \ddots & & \\
			& & & \ddots & & \\
			& & & & * & 0 \\
			& & & & -b_{d_b,d_b-1}&0 
		\end{array}
		\right].
	\end{equation}
	Laplace expansion about the final column gives
	\begin{align}
		\det\left(\mathbf{C}_{22}\xleftarrow[d_b]{}\mathbf{e}_1\right) &= (-1)^{d_b-1}(-1)^{d_b-1}\prod^{d_b-1}_{j=1}b_{j+1,j} \nonumber \\
		&= \prod^{d_b-1}_{j=1}b_{j+1,j}.
	\end{align}
	Therefore,
	\begin{equation}
		v = \cfrac{\prod^{d_b-1}_{j=1}b_{j+1,j}}{b(z)}
	\end{equation}
	because $\det \mathbf{C}_{22} = \det\left(z\mathbf{I} - \mathbf{B}\right) = b(z)$ by hypothesis.
	\item Now
	\begin{equation}
		\mathbf{C}_{12}\mathbf{C}_{22}^{-1}\mathbf{C}_{21} =
		\left[
		\begin{array}{cccc}
			& & & \alpha c_0 \\
			& & & \\
			& & & \\
			& & & \\
		\end{array}
		\right]
		\left[
		\begin{array}{cccc}
			& & & * \\
			& & & \vdots \\
			& & & * \\
			& & & v
		\end{array}
		\right] =
		\left[
		\begin{array}{cccc}
			& & & \alpha c_0v \\
			& & & \\
			& & & \\
			& & & 
		\end{array}
		\right]
	\end{equation}
	is $d_{a}+1$ by $d_a + 1$ and has its only nonzero entry, $\alpha c_0 v$, in the upper right corner.
	\item The Schur complement is therefore
	\begin{equation}
		\newcommand*{\tempb}{\multicolumn{1}{:c}{\begin{array}{c} -\alpha c_0 v \\ 0 \\ \vdots \\ 0 \end{array}}}
		\left[
		\begin{array}{cc}
			z\mathbf{I} - \mathbf{A} & \tempb \\ \cdashline{1-1}
			\begin{array}{cccc}
				0&\cdots&0&1
			\end{array}&
			z
		\end{array}
		\right]
	\end{equation}
	and we compute $\det\left(\mathbf{C}_{11} - \mathbf{C}_{12}\mathbf{C}_{22}^{-1}\mathbf{C}_{21}\right)$ by Laplace expansion on the last column:
	\begin{align}
		\det\left(\mathbf{C}_{11} - \mathbf{C}_{12}\mathbf{C}_{22}^{-1}\mathbf{C}_{21}\right) =& -(-1)^{d_a}\alpha c_0 v\det\left[
		\begin{array}{ccccc}
			-a_{21} & * & * & \cdot & * \\
			& -a_{32} & * & & *\\
			& & -a_{43} & & \vdots \\
			& & & \ddots & \\
			& & & & -a_{d_a,d_a-1} 
		\end{array}
		\right] \nonumber \\
		& + z\det\left(z\mathbf{I}-\mathbf{A}\right) \nonumber \\
		=& -(-1)^{d_a}\alpha c_0 v \prod^{d_a-1}_{j=1}\left(-a_{j+1, j}\right) + z\cdot a(z) \nonumber \\
		=& \alpha v \prod^{d_{a}-1}_{j=1} a_{j+1,j}\cdot c_0 + z\cdot a(z) \nonumber \\
		=& \alpha \cdot \cfrac{\left(\prod^{d_b-1}_{j=1}b_{j+1, j}\right)}{b(z)}\cdot\left(\prod^{d_a-1}_{j=1}a_{j+1, j}\right) \cdot c_0 + z\cdot a(z) \nonumber \\
		=& \cfrac{c_0}{b(z)} + z \cdot a(z)
	\end{align}
	by the definition of $\alpha$.
	
	Therefore by the Schur determinantal formula
	\begin{align}
		\det \left(z\mathbf{I}-\mathbf{C}\right) & = \det \left(\mathbf{C}_{22}\right)\det\left(\mathbf{C}_{11} - \mathbf{C}_{12}\mathbf{C}_{22}^{-1}\mathbf{C}_{21}\right) \nonumber \\
		&= b(z) \left(\frac{c_0}{b(z)} + z\cdot a(z)\right) \nonumber \\
		&= z \cdot a(z) b(z) + c_0.
	\end{align}
	Since the left hand side is a polynomial as is the right hand side, the formula will be true even if $b(z) = 0$, by continuity.
	\end{enumerate}
\end{proof}

If $p_n = \det\left(z\mathbf{I} - \mathbf{M}_n\right)$ for the Mandelbrot polynomials, the subdiagonals are all $-1$ and the matrices are the same size, so $\alpha = 1$ as is $c_0$: $p_{n+1} = zp_n^2 + 1$ gives
\begin{equation}
	\mathbf{M}_{n+1} = \left[
	\begin{array}{ccc}
	\mathbf{M}_n & & -\mathbf{c}_n\mathbf{r}_n \\
	-\mathbf{r}_n & 0 & \\
	& -\mathbf{c}_n & \mathbf{M}_n
\end{array}
	\right],
	\label{eq: Piers}
\end{equation}
where $\mathbf{r}_{n} = \left[ \begin{array}{cccc} 0 & 0 & \dots & 1 \end{array} \right]$ and $\mathbf{c}_n = \left[ \begin{array}{cccc} 1 & 0 & \cdots & 0\end{array}\right]^T$ are both of length $d_n$. This is Piers Lawrence's original construction \cite{piers}.These are remarkable matrices: they contain only $-1$ or $0$, and therefore are Bohemian matrices; yet the characteristic polynomial contains coefficients that grow exponentially in the degree $d_n$ (doubly exponentially in $n$).

For the Fibonacci-Mandelbrot polynomials, $\deg q_n = F_n - 1$ and the construction contains matrices of different size:
\begin{equation}
	\mathbf{M}_{n+1}=
  		\left[
  		\begin{array}{ccc}
  			\mathbf{M}_{n} & & (-1)^{d_{n+1}}\mathbf{c}_{n}\mathbf{r}_{n-1} \\
  			-\mathbf{r}_{n} & 0 & \\
  			& -\mathbf{c}_{n-1} & \mathbf{M}_{n-1}
  		\end{array}
  		\right],
\end{equation}
where $\mathbf{r}_{n} = \left[ \begin{array}{cccc} 0 & 0 & \cdots & 1 \end{array} \right]$ and $\mathbf{c}_{n} =  \left[ \begin{array}{cccc} 1 & 0 & \cdots & 0 \end{array} \right]^T$ be row and column vectors of length $d_{n}$. This gives a matrix of slightly greater height than \ref{eq: Piers} because entries may be ${-1, 0, 1}$.

For the Narayana-Mandelbrot polynomials, the product of $d_n-1$ $(-1)$s with $d_{n-2}-1$ $(-1)$s gives $(-1)^{d_n+d_{n-2}} = (-1)^{d_{n+1}}$ again.

This construction allows new matrix families. Suppose $s_0 = 0$, $s_{n+1} = z^3s_n^4 + 1$. Then if $\mathbf{S}_n$ is an upper Hessenberg companion for $s_n$ (with all $-1$ on the subdiagonal) the matrix
\begin{equation}
	\mathbf{S_{n+1}} = \left[
		\begin{array}{ccccccc}
			\mathbf{S}_n & & & & & & -\mathbf{c}_n\mathbf{r}_n \\
			-\mathbf{r}_n & 0 & & & & & \\
			& -\mathbf{c}_n & \mathbf{S}_n & & & & \\
			& & -\mathbf{r}_n & 0 & & & \\
			& & & -\mathbf{c}_n & \mathbf{S}_n & & \\
			& & & & -\mathbf{r}_n & 0 & \\
			& & & & & -\mathbf{c}_n & \mathbf{S}_n 
		\end{array}
	\right]
\end{equation}
is an upper Hessenberg companion for $s_{n+1}$. 

\section{Concluding Remarks}
This is a genuinely new kind of companion matrix. We demonstrate this on Newton's example polynomial $x^3 - 2x - 5$. We see that $x^3 - 2x - 5 = x(x^2 -2)-5 = x(x-\sqrt{2})(x+\sqrt{2})-5$, and companion matrices for $x-\sqrt{2}$ and $x + \sqrt{2}$ are just $[+\sqrt{2}]$ and $[-\sqrt{2}]$ respectively. Thus a companion matrix for Newton's polynomial is
\begin{equation}
	\left[
	\begin{array}{ccc}
		\sqrt{2} & & 5 \\
		-1 & & \\
		& -1 & -\sqrt{2}
	\end{array}
	\right]
\end{equation}
For unimodular polynomials, such companion matrices will be of lower height than the Frobenius or Fiedler \cite{fiedler2003note} companions, and may offer better numerical condition. 

We have now established that if $c(z) = z\cdot a(z)b(z) + c_0$ and $\mathbf{A}$ and $\mathbf{B}$ are upper Hessenberg companion matrices for the polynomials $a(z)$ and $b(z)$ respectively, then
\begin{equation}
	\mathbf{C} = \left[
	\begin{array}{ccc}
		\mathbf{A} & & -\alpha c_0 \mathbf{c_a}\mathbf{r_b} \\
		-\mathbf{r_a} & 0 & \\
		& -\mathbf{c_b} & \mathbf{B}
	\end{array}
	\right]
\end{equation}
is a companion matrix for $c(z)$. One wonders immediately about a corresponding linearization, $\mathbf{L_C}$, strong or otherwise, for the Matrix polynomial
\begin{equation}
	\mathbf{C}(z) = z \mathbf{A}(z) \mathbf{B}(z) + \mathbf{C_0}
\end{equation}
if $\mathbf{L_A}$ is a linearization for $\mathbf{A}$, $\mathbf{L_B}$ for $\mathbf{B}$. Some very preliminary experiments, where $\mathbf{L_A}$ and $\mathbf{L_B}$ were block upper Hessenberg with all blocks $\mathbf{I}$, so $\alpha = 1$, find that indeed
\begin{equation}
	\mathbf{L_C} = \left[
	\begin{array}{ccc}
	\mathbf{L_A} & & -c_0 \\
	\begin{array}{ccc}
	&&-\mathbf{I}
	\end{array} & 0 & \\
	& \begin{array}{c}
	-\mathbf{I} \\
	{ } 
	\end{array}
	& \mathbf{L_B}
	\end{array}
	\right]
\end{equation}
is a (strong) linearization for $c(z)$, in the examples we tried.

But we have no proof, and there are complications that suggest care will need to be taken. For instance, the matrix polynomials $\mathbf{C}_1 = z\mathbf{AB} + \mathbf{C}_0$ and $\mathbf{C}_2 = z\mathbf{BA}+\mathbf{C}_0$ may be different and have different polynomials eigenvalues. Placement of $\mathbf{L_B}$ in the lower right seems to be necessary, and different to exchange of $\mathbf{L_A}$ and $\mathbf{L_B}$.

We leave this extension to future work.

\section{Acknowledgements}
This work was supported by NSERC and by The University of Western Ontario, aka Western University. We thank Neil J. A. Sloane for introducing us to Narayana sequences, and Dario Bini \cite{bini2000design, bini2014solving} for teaching us about Mandelbrot polynomials and the Schur complement.

\clearpage
\section*{References}

\bibliography{mybibfile}

\end{document}